\documentclass{amsart}

\usepackage{amsmath,amsfonts,amssymb,amsthm}
\usepackage[mathscr]{eucal} 

\newtheorem{theorem}{Theorem}[section]
\newtheorem{lemma}[theorem]{Lemma}

\theoremstyle{definition}
\newtheorem{definition}[theorem]{Definition}

\newtheorem{example}[theorem]{Example}
\newtheorem{proposition}[theorem]{Proposition}

\theoremstyle{remark}

\newcommand{\seq}{$ (x_i)^\infty_{i=0}$} 

\begin{document}

\title{Transitivity and Mixing Properties of Set-valued Dynamical Systems}

\author{Wong Koon Sang}
\address{School of Informatics and Applied Mathematics,
	     Universiti Malaysia Terengganu, 21030 Kuala Nerus, Terengganu, Malaysia}
\email{kswong0910@gmail.com}

\author{Zabidin Salleh}
\address{School of Informatics and Applied Mathematics,
	     Universiti Malaysia Terengganu, 21030 Kuala Nerus, Terengganu, Malaysia}
\email{zabidin@umt.edu.my}





\begin{abstract}
We introduce and study two properties of dynamical systems: topologically transitive and topologically mixing under the set-valued setting. We prove some implications of these two topological properties for set-valued functions and generalize some results from single-valued case to set-valued case. We also show that both properties of set-valued dynamical systems are equivalence for any compact intervals.
\end{abstract}

\maketitle

\section{Introduction}

In dynamical systems, one of the most important research topics is to determine the chaotic behaviour of the system. Various definitions of chaos were introduced (see \cite{devaney89, li75, auslander80, akin10, schweizer94}) but up till now there is no universally accepted definition of chaos. Most of the definitions of chaos are closely related to transitivity of the dynamical systems. The concept of topologically transitive was started by Birkhoff \cite{birkhoff68} back in 1920. Dynamical systems with topologically transitive property contain at least one point which moves under iteration from one arbitrary neighborhood to any other neighborhood. Topologically transitive is one of the properties in topological dynamics that commonly use since it is a global characteristic in the dynamical system. Some prefer to study the topologically mixing of dynamical systems as it covers transitivity property as well.

Numerous studies related to the transitivity and mixing properties of the dynamical systems especially in one-dimensional have been done, see \cite{coven86, vellekoop94, blokh87, banks97, barge85, barge87}. As we all know, usually dynamical systems are studied in the view of single point. However, knowing how the points of the systems move is not sufficient as there are cases or problems that require one to know how the subsets of the system move. Loranty and Pawlak \cite{loranty12} studied the connection between transitivity and dense orbit for multifunction in generalized topological spaces. In recent years, several works and research on the topological dynamics of set-valued dynamical systems can be found (see \cite{woods18, liao06, roman03, bauer75, li15, metzger17, kennedy18, cordeiro16}). However, many properties for the dynamics of set-valued dynamical systems are yet to be discovered. More information about the transitivity in single-valued dynamical systems can be found in \cite{glasner06, du18, akin12, kolyada97}.

In this paper, we will introduce and study the notion of topologically transitive and topologically mixing for set-valued functions. We prove some elementary results of these two properties. Some of the results are generalization from the single-valued case (for e.g. \cite{ruette18, block92, coven86, holmgren94}). We also prove that the definitions of these two properties for set-valued function on compact intervals are equivalence. This paper is organized as follows. In Section 2 we give some background settings and define topologically transitive and topologically mixing for set-valued functions. In Section 3 we present some elementary implication results of topologically transitive and topologically mixing. In Section 4 we prove the equivalence of topologically transitive and topologically mixing of set-valued function on arbitrary compact interval. In Section 5 we present some conclusions.

\section{PRELIMINARIES}	\label{sec 2}

Let $X$ be a compact metric space. We denote $2^X$ as the collection of all nonempty closed subsets of $X$. We call a function $F: X \to 2^X$ as set-valued function. If $A\subseteq X$, then $F(A) = \{y \in X : \text{there is a point} \  x \in A \ \text{such that} \ y \in F(x)\}$. $F$ is said to be upper semicontinuous at $x \in X$ if for any open subset $V$ of $X$ containing $F(x)$, there is an open subset $U$ of $X$ containing $x$ such that for every $t \in U$, $F(t) \subseteq V$. $F$ is upper semicontinuous if it is upper semicontinuous at every point of $X$. Throughout the paper we assume the set-valued function $F$ is upper semicontinuous unless explicitly stated.

Since $X$ is compact, by \cite[(0.8)]{nadler78}  the hyperspace $2^X$ is compact. Therefore every element of $2^X$ is a nonempty compact subset of $X$ since they are non-empty closed subsets of $X$. The pair $(X, F)$ is called as set-valued dynamical system. $F^0$ is denoted as the identity on $X$ and $F^n = F \circ F^{n-1}$ for all integers $n > 0$.

Recall that in single-valued dynamical system $(X,f)$ where $f : X \to X$ represent a continuous function, for any point $x\in X$ we define the orbit of $x$ under $f$ as \seq \ where $x_0 = x$ and $x_{i+1} = f (x_i)$ for all integers $i \geq 0$. The point $x \in X$ is said to be a periodic point of $f$ with period $n$ provided $f^n (x) =x$ and $ f^j (x) \neq x$ for all integers $ 0 < j < n$. If the point $x$ has period $n = 1$ then it is called as fixed point. We extend these definitions to set-valued case.

\begin{definition} [\cite{raines18}] \label{orbit}
Let $(X, F)$ be a set-valued dynamical system. For any point $x \in X$, an orbit of $x$ is a sequence \seq \ such that $x_0 = x$ and $ x_{i+1} \in F (x_i)$ for all integers $ i \geq 0$. The collection of all orbits of $x$ is called as complete orbit of $x$, denoted by $CO (x)$.
\end{definition}

\begin{definition} [\cite{raines18, ansari10}] \label{p orbit}
For a set-valued dynamical system $(X, F)$ let  $x \in X$ and let \seq \ be an orbit of $x$. The orbit is said to be a periodic orbit if there exists $m \in \mathbb{N}$ such that $x_i = x_{i+m}$ for all integers $i \geq 0$. The point $x$ is a periodic point if it has at least one periodic orbit. The period of $x$ is the smallest number $m \in \mathbb{N}$ satisfying $x_i = x_{i+m}$ for all integers $ i \geq 0$. If $m =1$ then $x$ is said to be a fixed point.
\end{definition}

From Definitions \ref{orbit} and \ref{p orbit}, we can see that in set-valued dynamical systems, the orbits of $x$ under $F$ no longer uniquely determined. The following example shows that the orbit \seq \ is not necessarily periodic even if there exists $j \in \mathbb{N}$ such that $x= x_0 = x_j$.

\begin{example}
Let $X = \{0, 1\}$ and let the set-valued function $F: X \to 2^X$ defined by $F (1) = \{ 0,1 \}$ and $F (0) = \{ 1 \}$. Let $x=0$, then one of the orbit of $x$ under $F$ is $ (0, 1, 1, 0, 1, 0, 1,  \ldots )$. We can see that $x =x_3 $ but $ x_i \neq x_{i+3}$ for some $i > 0$.
\end{example}

Next we define topologically transitive of set-valued functions. For topologically mixing of set-valued functions we adopt the definition which has been defined by \cite{raines18}. Note that the set-valued product function $F \times F : X \times X \to 2^{X \times X}$ is defined by $ \left(F \times F\right) ( x, x') = \{ (y, y') \in X \times X : y \in F(x) \text{ and } y' \in F (x') \}$.

\begin{definition} \label{trans}
A set-valued function $F$ is topologically transitive if for any nonempty open subsets $U$ and $V$ of $X$, there exists $m \in \mathbb{N}$ and $x \in U$ with an orbit \seq \ such that $x_m \in V$.
\end{definition}

\begin{definition} [\cite{raines18}] \label{mix}
A set-valued function $F$ is topologically mixing if for any nonempty open sets $U$ and $V$ in $X$, there is an $M \in \mathbb{N}$ such that for any $m > M$ there is an $x_0 \in U$ with an orbit \seq \ such that $x_m \in V$.
\end{definition}

\begin{definition} \label{vtrans}
Let $F$ be a set-valued function of $X$. Then $F$ is said to be topologically bitransitive if $F^2$ is topologically transitive. $F$ is totally transitive if $F^n$ is topologically transitive for all $n \in \mathbb{N}$. $F$ is topologically weakly mixing if the set-valued product function $F \times F$ is topologically transitive.
\end{definition}

We end this section by recall some concepts from topology. Let $U$ be any subset of $X$, then the interior and the closure of $U$ is denoted by $\text{int} ( U)$ and $\text{cl}( U)$ respectively. A set $U$ is said to be dense in $X$ if $\text{cl}( U) = X$. In other words, we can say that $U$ is dense in $X$ if every open subset of $X$ contains at least a point of $U$ (see \cite{munkres00} and \cite{roseman99}).

\section{TOPOLOGICALLY TRANSITIVE AND MIXING OF SET-VALUED FUNCTION} \label{sec 3}

In single-valued case, there are two commonly used definitions for topologically transitive: one is defined by using open sets and another one is defined by using points with dense orbit (see \cite{kolyada97}). Block \cite{block92} showed that both definitions coincide when the space is compact. But in general, both characterizations of transitivity are not equivalent as shown in \cite{peris99} and \cite{degirmenci03}. On a compact metric space with set-valued function, we show that if there is a point with dense orbit then it will imply the transitivity of the set-valued function.

\begin{proposition} \label{point dense}
Let $F: X \to 2^X$ be a set-valued function. If there exists at least a point $x \in X$ with an orbit  \seq \ $\in CO(x)$ such that the orbit is dense in $X$, then $F$ is topologically transitive.
\end{proposition}
\begin{proof}
Let $U$ and $V$ be open sets in $X$. We need to find a point $x \in U$ and a positive integer $m$ such that $x$ has an orbit \seq \ with $x_m \in V$. By hypothesis, there is a point in $X$ with an orbit that is dense in $X$. Let $y \in X$ be such a point with an orbit $ (y_i)^\infty_{i=0}$, then we know that there exists a positive integer $k$ such that $y_k \in U$. Now we try to show that there is a positive integer $m$ such that $y_{k+m} \in V$. Then the proof is done by letting $x = y_k$, so we have $x \in U$ with an orbit \seq \ where $x_m = y_{k+m} \in V$.

Since the orbit $ (y_i)^\infty_{i=0}$ of $y$ is dense, $V$ contains at least one iterate of $y$. Suppose there are only finitely many iterates of $y$ in $V$. Let $z$ be any element of $V$ such that $z$ is not an iterate of $y$ and let $\varepsilon = \text{min}\{ d(z,y_j) : j = 0, 1, 2, 3, \ldots \}$ where $d$ is the metric on $X$. We have $\varepsilon > 0$ and the neighborhood $N_{\varepsilon /2} ( z)$ does not contain any iterates of $y$. This implies that the orbit $(y_i)^\infty_{i=0}$ of $y$ is not dense in $X$, a contradiction. Therefore, $V$ must contain infinitely many iterates of $y$. Since there are only finitely many positive integers less than $k$, there exists an integer $n > k$ such that $y_n \in V$. We let $m = n - k$ and the proof is complete.
\end{proof}

With Proposition \ref{point dense}, we obtain the following theorem.

\begin{theorem} \label{equivalent}
A set-valued function $F: X \to 2^X$ is topologically transitive if and only if there exists at least a point $x \in X$ with an orbit \seq $\in CO(x)$ such that the orbit is dense in $X$.
\end{theorem}
\begin{proof}
Clearly by the definition of topologically transitive, we will have at least a point $x \in X$ with an orbit  \seq $\in CO(x)$ such that the orbit is dense in $X$. For the converse part we have proved in Proposition \ref{point dense}.
\end{proof}

Similar to the single-valued case (see \cite{denker76}), it is easy to see that by the definition, if the set-valued function $F$ is topologically mixing then it implies that $F$ is topologically transitive. We show that the converse is not true in the following example.

\begin{example}
Let us consider the unit circle $S^1 = \{ z \in \mathbb{C} : | z | = 1 \}$ and $R_\alpha : S^1 \to S^1$  be an irrational rotation of $S^1$. We define a set-valued function $R$ of $S^1$ by $R ( z) = [ z, R_\alpha(z)]$. Since all the points of $S^1$ have dense orbits, for an open subset $U$ of $S^1$, it will intersect with other open subset $V$ of $S^1$ for some iterates under $R$. Therefore $R$ is topologically transitive. But $R_\alpha$ is irrational rotation, so there exists at least one further iterates of $U$ under $R$ that did not intersect with $V$. Hence $R$ is not topologically mixing.
\end{example}

Next we discuss some connections between Definition \ref{trans}, \ref{mix} and \ref{vtrans} in Section \ref{sec 2}.

\begin{proposition}
Let $(X,F)$ be a set-valued dynamical system. If the set-valued function $F$ is topologically mixing, then $F$ is topologically weakly mixing.
\end{proposition}
\begin{proof}
Assume that $F$ is topologically mixing. Let $ W_1, W_2$ be any two nonempty open sets in $X \times X$. There exist nonempty open sets $U_1, U_2, V_1, \text{and } V_2$ in $X$ such that $U_1 \times U_2 \subset W_1$ and $V_1 \times V_2 \subset W_2$. Since $F$ is topologically mixing, there exists $N_1 \in \mathbb{N}$ such that for all $n \geq N_1$, there is an $x \in U_1$ with an orbit \seq \ such that $x_n \in V_1$. Similarly, there exists $N_2 \in \mathbb{N}$ such that for all $n \geq N_2$, there is an $y \in U_2$ with an orbit $ (y_i)^\infty_{i=0}$ such that $y_n \in V_2$. Let $M = \text{max} \{ N_1, N_2 \}$. Then for all $n \geq M$, there exists $ (x , y) \in W_1$ with an orbit $ (x_i,y_i)^\infty_{i=0}$ such that $ (x_n , y_n ) \in W_2$. Hence we conclude that $F$ is topologically weakly mixing.
\end{proof}

\begin{proposition}
Let $(X,F)$ be a set-valued dynamical system. If the set-valued function $F$ is topologically mixing, then $F$ is topologically bitransitive.
\end{proposition}
\begin{proof}
Assume that $F$ is topologically mixing. Then for any two nonempty open sets $U$ and $V$ of $X$, there is an $M \in \mathbb{N}$ such that for any positive integer $m > M$, there is an $x \in U$ with an orbit \seq \ such that $x_m \in V$. If we take $m$ to be an even positive integer greater than $M$, i.e. $m =2n > M$ where $n$ is a positive integer, then there exists an $x \in U$ with an orbit \seq \ such that $x_{2n} \in V$. This implies that $F^2$ is topologically transitive and hence $F$ is topologically bitransitive.
\end{proof}

\begin{lemma} \label{product trans}
Let $(X,F)$ be a set-valued dynamical system and the intersection of an open set with any iterates of another open set contains an open set of $X$. If the set-valued function $F$ is topologically weakly mixing, then the set-valued product dynamical system $( X^n, F \times \cdots \times F)$ is topologically transitive for all integers $n \geq 1$.
\end{lemma}
\begin{proof}
For all open sets $U$, $V$ in $X$, we define
$$ N ( U, V) = \{ n\geq 1 : \text{ there exists } x \in U \text{ with an orbit } (x_i)^\infty_{i=0}  \text{ such that } x_n \in V \}. $$
Let $U_1, U_2, V_1$ and $V_2$ be nonempty open sets in $X$. Since $F$ is topologically weakly mixing, there exists a natural number $m$ such that there is a point $ ( x, y ) \in U_1 \times V_1$ with an orbit $( x_i, y_i )^\infty_{i=0}$ such that $( x_m, y_m ) \in U_2 \times V_2$. This means that there is a point $ x \in U_1$ with orbit $ ( x_i)^\infty_{i=0}$ such that $x_m \in U_2$ and there is a point $ y \in V_1$ with orbit $ (y_i)^\infty_{i=0}$ such that $y_m \in V_2$. We can say that for any $G, H$ nonempty open sets in $X$,  $N ( G, H) \neq \emptyset$.

Now we are going to show that there exist nonempty open sets $U, V$ in $X$ such that $N( U, V) \subseteq N (U_1, V_1) \cap N (U_2, V_2)$. Let us define the open sets $U \subseteq U_1 \cap F^n(U_2)$ and $V \subseteq V_1 \cap F^n (V_2)$ as follow:
$$U := \{ x : x \in U_1 \text{ with orbit } (x_i)^\infty_{i=0} \text{ such that } x_0 =x \text{ and } x_n \in U_2 \}$$
and
$$V := \{ x : x \in V_1 \text{ with orbit } (x_i)^\infty_{i=0} \text{ such that } x_0 =x \text{ and } x_n \in V_2 \}.$$
We have already shown that these sets are not empty. Let $k \in N (U, V)$. This integer exists and satisfy the condition of there exists an $x \in U$ with an orbit \seq \ such that $x_0 =x$ and $x_k \in V$. This means that there is an $x \in U$ with an orbit \seq \ such that $x_0 =x \in U_1$, $x_n \in U_2, x_k \in V_1 \text{ and } x_{k+n} \in V_2$. Then we can deduce that there is an $x \in U_1$ with orbit \seq \ such that $x_k \in V_1$ and there is an $x \in U_2$ with an orbit \seq \ such that $x_k \in V_2$. Therefore, we obtain $k \in N (U_1, V_1)$ and $k \in N(U_2, V_2)$ and this implies that $N( U, V) \subseteq N (U_1, V_1) \cap N (U_2, V_2)$. By using principal of mathematical induction, we able to see that for all nonempty open sets $U_1, \ldots, U_n, V_1, \ldots, V_n$ in $X$, there exist nonempty open sets $U, V$ in $X$ such that
$$N( U, V) \subseteq N (U_1, V_1) \cap N (U_2, V_2) \cap \cdots \cap N (U_n, V_n).$$
Hence, we conclude that $( X^n, F \times \cdots \times F)$ is topologically transitive.
\end{proof}

\begin{theorem}
Let $(X,F)$ be a set-valued dynamical system and the intersection of an open set with any iterates of another open set contains an open set of $X$. If the set-valued function $F$ is topologically weakly mixing, then $F$ is totally transitive.
\end{theorem}
\begin{proof}
Assume that $F$ is topologically weakly mixing. Let $n \geq 1$ be a fixed positive integer and $U, U', V, V'$ be nonempty open sets in $X$. We define two open sets $W, W'$ in $X^{2n}$ as follow:
$$ \begin{aligned}
W =& \{ (x_0, x_1, \ldots, x_{n-1}, y_0, y_1, \ldots, y_{n-1} ) : x_i \in F^i (U) , \\
& y_i \in F^i (V) \text{ for all } i=0, 1, \ldots, n-1 \}
\end{aligned}$$
and
$$ \begin{aligned}
W' =& \{ (x'_0, x'_1, \ldots, x'_{n-1}, y'_0, y'_1, \ldots, y'_{n-1} ) : x'_i \in U' , \\
& y'_i \in V' \text{ for all } i=0, 1, \ldots, n-1 \}.
\end{aligned}$$
By Lemma \ref{product trans}, $(X^{2n}, F\times F \times \cdots \times F)$ is topologically transitive. So, there exist a positive integer $m$ such that there is an $(x_0, x_1, \ldots, x_{n-1}, y_0, y_1, \ldots, y_{n-1} ) \in W$ with an orbit $\left (x^{(i)}_0, x^{(i)}_1, \ldots, x^{(i)}_{n-1}, y^{(i)}_0, y^{(i)}_1, \ldots, y^{(i)}_{n-1} \right ) ^\infty_{i=0}$ such that
$$\left (x^{(0)}_0, x^{(0)}_1, \ldots, x^{(0)}_{n-1}, y^{(0)}_0, y^{(0)}_1, \ldots, y^{(0)}_{n-1} \right ) \in W$$
and
$$\left (x^{(m)}_0, x^{(m)}_1,  \ldots, x^{(m)}_{n-1}, y^{(m)}_0, y^{(m)}_1, \ldots, y^{(m)}_{n-1} \right ) \in W'.$$
 This implies that for all $i \in \{0, 1 , \ldots , n-1 \}$, there is an $x \in U$ with an orbit \seq \ such that $x_0 \in U$ and $x_{i+m} \in U'$ and there is an $y \in V$ with an orbit $ (y_i)^\infty_{i=0}$ such that $y_0\in U$ and $y_{i+m} \in V'$. We choose an $i \in \{0, 1 , \ldots , n-1 \}$ such that $m+i = np$ where $p$ is a positive integer. We can conclude that there is a point $(x,y) \in U \times V$ with an orbit $ (x_i, y_i)^\infty_{i=0}$ such that $(x_{np}, y_{np}) \in U' \times V'$. Therefore, $F^n$ is topologically weakly mixing which implies that $F^n$ is topologically transitive and completes the proof.
\end{proof}

\begin{proposition}
Let $(X,F)$ be a set-valued dynamical system. If the set-valued function $F$ is totally transitive, then $F$ is topologically bitransitive.
\end{proposition}
\begin{proof}
Let $U,V$ be any two nonempty open sets of $X$. Since $F$ is totally transitive, $F^n$ is topologically transitive for all positive integers $n \geq 1$. For each $n \geq 1$, there exists an $m \in \mathbb{N}$ such that there is an $x\in U$ with an orbit \seq \ such that $x_{mn} \in V$. When we take $n=2$, there exists an $m \in \mathbb{N}$ such that there is an $x\in U$ with an orbit \seq \ such that $x_{2m} \in V$. This implies that $F^2$ is topologically transitive, hence $F$ is topologically bitransitive.
\end{proof}

The implications between the various conditions of topologically transitive and mixing in this section are summarized as follows:
\bigskip
$$\text{mixing} \Longrightarrow \text{weakly mixing} \Longrightarrow \text{totaly transitive} \Longrightarrow \text{bitransitive} \Longrightarrow \text{transitive}$$

\bigskip
\section{TRANSITIVITY AND MIXING PROPERTIES OF SET-VALUED FUNCTION IN COMPACT INTERVALS}

In this section, we investigate the properties of topologically transitive and topologically mixing of set-valued functions on arbitrary compact interval $I=[a,b]$, where $a, b \in \mathbb{R}$ such that $a < b$. In single-valued case,  the set of periodic points is dense in $I$ if the function is topologically transitive \cite{block92}. By the help of the following lemma, similar result can be obtained in set-valued case.

\begin{lemma} \label{inequality of point}
Let $(I,F)$ be a set-valued dynamical system and let $J$ be a subinterval of $I$ which contains no periodic point of $F$. Suppose that , $x \in J$ has an orbit \seq \ such that $x_m \in J$ for some integers $m > 0$ and $y \in J$ has an orbit $(y_i)^\infty_{i=0}$ such that $y_n \in J$ for some integers $n > 0$. If $x < x_m$ then $y < y_n$ and if $x > x_m$ then $y > y_n$.
\end{lemma}
\begin{proof}
Without loss of generality, suppose that $x < x_m$. Let $G = F^m$, then the subinterval $[x, x_m]$ of $J$ contains no periodic point of $G$. If $x_{km} >x$ for some $k \geq 1$ then $x_{(k+1)m} > x_m$ as there is no fixed point of $G^k$ inside the interval $[x, x_m]$. Clearly by mathematical induction, $x_{km} > x$ for all $k \geq 1$. So, in particular we have $ x_{mn} > x$.

Assume that $y > y_n$. With similar argument as above (the order of inequality is reversed) we yield $y >y_{mn}$. Consequently, there exists a point $z \in J$ lies in between $x$ and $y$ such that $z_{mn} = z$. This leads to a contradiction as $J$ contains no periodic point of $F$. Therefore, we conclude that $y < y_n$.
\end{proof}

\begin{proposition}
If the set-valued function $F:I \to 2^I$ is topologically transitive, then the set of periodic points of $F$ is dense in $I$.
\end{proposition}
\begin{proof}
We will prove this by contradiction. Suppose that the set of periodic points of $F$ is not dense in $I$. There exist two points $x, y \in I$ where $ x < y$ such that the open interval $ (x, y)$ contain no any periodic points of $F$. Since $F$ is topologically transitive, by Theorem \ref{equivalent} there is a point $ u \in (x, y)$ with a dense orbit $(u_i)^\infty_{i=0}$ in $I$. Hence, for some integers $m > 0$ and $0 < p < q$, we have $ u < u_m < y$ and $x < u_q < u_p < u$. Let $v=u_p$. Since $u_q \in F^{q-p} (v)$, $v$ has an orbit $(v_i)^\infty_{i=0}$ such that $v_{q-p} = u_q$. Then, we obtain the inequality
$$ x < v_{q-p} < v < u < u_m < y.$$
But this is impossible as Lemma \ref{inequality of point} applied to the open interval $(x, y)$. Thus we conclude that the set of periodic points of $F$ is dense in $I$.
\end{proof}

Next we prove an elementary property of topologically mixing for set-valued function in compact interval.

\begin{proposition} \label{mixing interval}
Let $I=[a, b]$ and $F:I \to 2^I$ be a set-valued function. $F$ is topologically mixing if and only if for all nondegenerate subintervals $J \subseteq I$ and any pair $c, d \in \text{int} (I)$ such that $ a < c < d < b$, there exists a positive integer $M$ such that $[c,d] \subset F^n (J)$ for all $n \geq M$.
\end{proposition}
\begin{proof}
Suppose that $F$ is topologically mixing. Let $\varepsilon >0$. Let $U = (a, c)$ and $V = (d, b)$. If $J$ is a nonempty open subinterval of $I$, then there exists $N_1 \in \mathbb{N}$ such that there is an $x \in J$ with an orbit \seq \ such that $x_n \in U$ for all $n \geq N_1$ since $F$ is topologically mixing. Similarly, there exists $N_2 \in \mathbb{N}$ such that there is an $x \in J$ with an orbit \seq \ such that $x_n \in V$ for all $n \geq N_2$. Let $M= \text{max} \{N_1, N_2 \}$. Then, we have an $x \in J$ with an orbit \seq \ such that $x_n \in U\cap V$ for all $n \geq M$. This implies that $[c, d] \subset F^n (J )$ by connectedness. If $J$ is a nondegenerate subinterval, the same result holds by consider the nonempty open interval $\text{int} (J)$.

Conversely, suppose that for any pair $c, d \in \text{int} (I)$ and a nondegenerate subinterval $J \subseteq I$, there is a positive integer $M$ such that $[c,d] \subset F^n (J)$ for all $n \geq M$. Let $U, V$ be two nonempty open subsets of $I$. Choose two nonempty open subintervals $J, K$ such that $J \subset U, K \subset V$ and neither $a$ nor $b$ is an endpoint of $K$. There exists a pair $c, d \in \text{int} ( I )$ such that $K \subset [c, d]$. By assumption, there exists positive integer $M$ such that $[c, d] \subset F^n (J)$ for all $n \geq M$. Therefore, we have $K \subset [c, d ]\subset F^n(J)$ and this implies that there is an $x \in U$ with an orbit \seq \ such that $x_n \in V$ for all $n \geq M$. We conclude that $F$ is topologically mixing.
\end{proof}

When topologically mixing is replace with topologically bitransitive for set-valued functions, a weaker version of Proposition \ref{mixing interval} can be obtain with the help of the following lemma. We refer the notation from \cite{ruette18}, $ \langle a, b \rangle$ is use to denote an interval with the pair of points $a$, $b$ as the endpoints and either $ a < b$ or $ a >b$.

\begin{lemma} \label{2 ways}
Let $F:I \to 2^I$ be a set-valued function. Let $J$ be a subinterval of $I$ which contains a fixed point $p$ of $F$ and a periodic point $q$ of $F$ with least period $m \geq 3$. Let $q_i$ represents the iterates of $q$ in the periodic orbit  where $i \in \{ 0, 1, \ldots , m-1 \}$. Then one of the following holds:
\begin{enumerate}
	\item $F^n (J) \supset \{q, q_1, \ldots, q_{m-1} \}$ for all $n \geq 2m-2$;
	\item $m$ is even, the set $\{q, q_2, \ldots, q_{m-2} \}$ and $\{q_1, q_3, \ldots, q_{m-1} \}$ lie on opposite sides of $u$ and $F^{2n} (J) \subset \{q, q_2, \ldots, q_{m-2} \}$ for all $n \geq (m-2)-1$.
\end{enumerate}
\end{lemma}
\begin{proof}
Let $u$ and $v$ be two adjacent points in the periodic orbits of $q$ such that $v' \leq u < v \leq u'$ where $v' \in F (v), u' \in F(u)$ and both $u', v'$ are contain in the same periodic orbit of $q$. Let $i, j \in \{ 0,1, \ldots, m-1 \}$ such that $u =q_i$ and $v = q_j$. If both point $q_k$ and $q_{k+1}$ lie on the same side of the fixed point $p$ for some $k \in \{ 0,1, \ldots, m-1 \}$, then we have four possible cases:
\renewcommand{\labelenumi}{(\roman{enumi})}
\begin{enumerate}
	\item $q_i = u < v < p$;
	\item $q_k < u < p < v$;
	\item $u < p < v < q_k$;
	\item $p < u < v = q_j$.
\end{enumerate} 
For Case (i), we have $F^i (J) \supset [u, v]$. For Case (iv), we have $F^i (J) \supset [u, v]$. For Case (ii) and (iii), the set $F^k(J)$ contains the compact interval $ \langle u, q_k \rangle$. Thus we have $F^{k+1} (J) = F \left ( F^k (J) \right ) \supset [u, v]$. In both cases, since $F( [u, v] ) \supset [u, v]$, we obtain $F^m( J ) \supset [u, v]$ and, since  $F^{m-2}( [u, v] ) \supset \{q, q_1, \ldots, q_{m-1} \}$, we have $F^n (J) \supset \{q, q_1, \ldots, q_{m-1} \}$ for all $n \geq 2m-2$.

Otherwise, if $q_i$ and $q_{i+1}$ lie on the opposite sides of $p$ for all $i \in \{ 0,1, \ldots, m-1 \}$, clearly $m$ is an even integer. Let $r$ be the even integer from the set $\{ 0,1, \ldots, m-1 \}$ such that $q_r$ is the point which stay on the same side of $z$ with $q$ and lie most far from $p$. Then, $F^{r+2j} (J)$ contains the set $ \{q_m =q, q_2,q_4, \ldots,  q_{m-2} \}$ for all $j \geq 0$. Therefore, we conclude that $F^{2n} (J) \supset \{q, q_2, \ldots, q_{m-2} \}$ for all $n \geq (m/2)-1$.
\end{proof}

\begin{theorem} \label{bitrans interval}
Let $I=[a, b]$ and $F: I \to 2^I$ be a set-valued function. If $F$ is topologically bitransitive then for any nondegenerate subinterval $J \subseteq I$ and any pair $c, d \in \text{int} (I)$ such that $a < c < d < b$, there exists a positive integer $M$ such that $[c, d] \subset F^n (J)$ for all $n \geq M$.
\end{theorem}
\begin{proof}
Let $z$ be a fixed point of $F$ and $J$ be any nondegenerate subinterval of $I$. Without the loss of generality, we may assume that $z <x$ for all $x \in J$. Since $F$ is topologically bitransitive, it means that $F^2$ is topologically transitive and by Theorem \ref{equivalent}, there exists a point $u \in \text{int} (J)$ with an orbit $(u_i)^\infty_{i=0}$ with respect to $F^2$ which is dense in $I$.

Let $L = [c, d]$ be a compact subinterval in $I$ where $c, d \in \text{int} (I)$ and $c < d$. Since the orbit $(u_i)^\infty_{i=0}$ with respect to $F^2$ which is dense in $I$, for some positive integers $n_1 < n_2$ we obtain
$$ u_{n_1} < \text{min}\{ z, \text{min } L \} <  \text{max}\{ z, \text{max } L \}< u_{n_2}.$$

By Proposition \ref{point dense}, we know that the set of periodic points is dense in $I$. Since $u \in \text{int} (K)$, there exists a periodic point $p \in K$ with period $m$ which close to point $u$ and the set contain all points of the periodic orbit of $p$ is contained in $\text{int}(I)$. Then, we have
$$ u_{n_1} \approx p_{n_1} < \text{min}\{ z, \text{min } L \} <  \text{max}\{ z, \text{max } L \}< p_{n_2} \approx u_{n_2}. $$

Consequently, there is a positive integer $s$ where the interval $F^s(K)$ contains the fixed point $z$ and the periodic point $p_s$ with period $m$. For the orbit $(p_i)^\infty_{i=0}$, we can see that the even iterates are distributed on both sides of $z$. Hence, by Lemma \ref{2 ways}, we have $F^j ( F^s (K) ) \supset \{ p, p_1, \ldots, p_{m-1} \}$ for all $j \geq 2m-2$. Therefore, $F^n (K) \supset [ \text{min} \{ p, p_1, \ldots, p_{m-1} \}, \text{max} \{ p, p_1, \ldots, p_{m-1} \} ] \supset L$ for all $ n \geq s+2m-2$ and the prof is complete.
\end{proof}

The following theorem gives an overview on the relation between topologically transitive and topologically mixing for set-valued functions of compact intervals. In fact, in compact intervals both definitions are equivalent, which similar to the results in single-valued case (see \cite{block92}).

\begin{theorem}
Let $I= [a, b]$ and $F: I \to 2^I$ be a set-valued function. If $F$ is topologically transitive, then the following statements are equivalent:
\begin{enumerate}
	\item $F$ is topologically bitransitive. \label{item1}
	\item $F$ is totally transitive. \label{item2}
	\item $F$ is topologically weakly mixing. \label{item3}
	\item $F$  is topologically mixing. \label{item4}
	\item For any nondegenerate subinterval $J \subseteq I$ and any pair $c, d \in \text{int}(J)$ such that $ a < c < d <b$, there exists a positive integer $M$ such that $[c ,d ] \subset F^n (J)$ for all $n \geq M$. \label{item5}
\end{enumerate}
\end{theorem}
\begin{proof}
It follows from Theorem \ref{bitrans interval} that $(1) \Rightarrow(5)$. By Proposition \ref{mixing interval}, we have $(5) \Leftrightarrow (4)$. Finally, by the implication diagram at the end of Section \ref{sec 3} we have $(4) \Rightarrow (3) \Rightarrow (2) \Rightarrow (1)$.
\end{proof}

\section{CONCLUSION}

In this paper we studied the two topological properties of  dynamical systems which are topologically transitive and topologically mixing under the setting of set-valued case. An implication diagram to show the connection between various conditions of transitivity and mixing is provided in Section \ref{sec 3}. We also investigated transitivity and mixing properties of set-valued functions for compact intervals and showed both definitions are equivalent, which similar to the single-valued case.

\bibliographystyle{plain}

\end{document}